\documentclass{article}

\usepackage{amsmath,amsthm,amssymb}
\usepackage{graphicx}

\newtheorem{theorem}{Theorem}

\newtheorem{lemma}[theorem]{Lemma}

\newtheorem{corollary}[theorem]{Corollary}

\newcommand{\gb}{Gr\"obner basis }
\newcommand{\gbs}{Gr\"obner bases }
\newcommand{\pe}{\preceq}
\newcommand{\se}{\succeq}

\def \cocoa{{\hbox{\rm C\kern-.13em o\kern-.07em C\kern-.13em o\kern-.15em A}} }

\begin{document}

\title{Computing Gr\"obner Bases of Ideals of Few Points in High Dimensions}

\author{Winfried Just$^{1}$ and Brandilyn Stigler$^2$\\ \\
{\small $^1$Department of Mathematics, Ohio University, Athens, OH 45701}\\
{\small $^2$Mathematical Biosciences Institute, The Ohio State
University, Columbus, OH 43210}}

\date{November 2, 2006}

\maketitle

\abstract{A contemporary and exciting application of \gbs is their
use in computational biology, particularly in the reverse
engineering of gene regulatory networks from experimental data.  In
this setting, the data are typically limited to tens of points,
while the number of genes or variables is potentially in the
thousands. As such data sets vastly underdetermine the biological
network, many models may fit the same data and reverse engineering
programs often require the use of methods for choosing parsimonious
models. \gbs have recently been employed as a selection tool for
polynomial dynamical systems that are characterized by maps in a
vector space over a finite field.

While there are numerous existing algorithms to compute Gr\"obner
bases, to date none has been specifically designed to cope with
large numbers of variables and few distinct data points.  In this
paper, we present an algorithm for computing \gbs of
zero-dimensional ideals that is optimized for the case when the
number~$m$ of points is much smaller than the number~$n$ of
indeterminates.  The algorithm identifies those variables that are
\emph{essential}, that is, in the support of the standard monomials
associated to a polynomial ideal, and computes the relations in the
\gb in terms of these variables. When $n$ is much larger than $m$,
the complexity is dominated by~$nm^3$. The algorithm has been
implemented and tested in the computer algebra system
\emph{Macaulay~2}. We provide a comparison of its performance to the
Buchberger-M\"oller algorithm, as built into the system.}


\section{Introduction}

The theory of \gbs has been an active field of study in the last
four decades, beginning with the seminal work of Buchberger
\cite{buchberger-phd}. A problem of particular interest has been the
development of algorithms for computing Gr\"obner bases. The first
algorithm, proposed by Buchberger, has time complexity that is
doubly exponential in the number of variables \cite{buchberger}.
Since then, several improvements to Buchberger's algorithm have been
proposed, as well as a number of alternative methods for certain
classes of ideals.

Many of the improvements focus on two aspects.  The first is
coefficient growth when computing \gbs in a field of characteristic
0 (for example, see \cite{arnold}). The second is Buchberger's
Criterion, which states that
\begin{quote}``A set $G=\{g_1,\ldots,g_r\}\subset I$ is a \gb for
$I$ if and only if the $S$-polynomial $\overline{S(g_i,g_j)}^G$ is 0
for all $1\leq i,j\leq r$.''
\end{quote}
The Optimized Buchberger Algorithm \cite{caboara} proposed by
Caboara \textit{et al.} and Faug\'ere's \emph{F}4 and \emph{F}5
\cite{faugere-f4,faugere-f5} are instances of methods that seek to
minimize the number of $S$-polynomials to be computed. While they
still have exponential complexity in the worst case, in practice
their performance renders them efficient alternatives to the
original Buchberger algorithm.

For zero-dimensional ideals, several methods have been described and
implemented.  In \cite{BM}, the authors presented the
Buchberger-M\"oller algorithm (BMA) for computing the reduced \gb
for the vanishing ideal of a variety $V$ over a field.  This
algorithm eliminates the need to compute $S$-polynomials and instead
performs Gaussian elimination on a generalized Vandermonde matrix.
Its complexity is quadratic in the number of variables and cubic in
the number of points in $V$ (for details, see \cite{mariani, mora,
robbiano}). It has been implemented in publicly available computer
algebra systems such as \cocoa \cite{cocoa} and \emph{Macaulay 2}
\cite{M2}. The BMA was later generalized to noncommutative rings
\cite{alonso}. Abbott \textit{et al.} \cite{abbott} described a
modular version of the BMA for the case when~$k=\mathbb{Q}$.

There are other algorithms for zero-dimensional ideals which have
been developed for particular settings. Farr and Gao presented an
algorithm based on a generalization of Newton interpolation in
\cite{farr}. While the complexity is exponential in the number $n$
of variables, the algorithm has been optimized for the case in which
$n$ is small as compared to the number of points. Lederer proposed a
method for lexicographic term orders which gives insight into the
structure of the Gr\"obner basis \cite{lederer}.

A recent and exciting development in the theory of \gbs is their use
in computational biology.  For instance they have been used in the
identification of critical points of maximum likelihood functions in
phylogenetic-tree reconstruction \cite{sturmfels}. \gbs have also
been employed as a selection tool for polynomial dynamical systems
(PDSs) in the study of gene regulatory networks \cite{LS} and
protein signal transduction networks~\cite{allen}.

In applications to molecular biology, networks often consist of $n$
biochemicals, such as gene products or metabolites, with changing
concentration levels.  In \cite{LS} a method was proposed to reverse
engineer biochemical networks, where the levels are mapped to a
finite field $k=\mathbb{F}_p$ for some prime $p$. In this setting,
networks are modeled as PDSs, which generalize the widely studied
Boolean networks (see \cite{kauffman} for an introduction).
Concentration levels are recorded in a vector in $k^n$, and the data
consists of input-output pairs $(s_i,s_{i+1})\in k^n\times k^n$,
where $s_i$ is a vector describing the state of the network at time
$i$, for $i=1,\ldots ,m$. The input vectors can be viewed as an
affine variety $V\subset k^n$, and a family of models represented as
PDSs is constructed in terms of the vanishing ideal of $V.$ \gbs are
then used to select the most parsimonious PDS from this collection.
In these applications, the number $n$ is typically in the hundreds
to thousands, whereas the number $m$ is at best on the order of tens
of measurements.

Below we describe an algorithm for computing \gbs for
zero-dimensional ideals (\textit{i.e.}, vanishing ideals) in a
polynomial ring $R$.  This algorithm is specialized for the case
when the number $m$ of distinct points is much smaller than the
number $n$ of variables. In this setting, there are few relations in
terms of \emph{essential} variables, that is, variables that are in
the support of the standard monomials associated to an ideal. The
remaining ones are of the type $x_i-g$ where the leading term $x_i$
is not an essential variable and the support of $g$ has only
essential variables. Therefore computation of a \gb can be
restricted to a proper subring of $R$ containing only essential
variables. The algorithm identifies these variables and computes
relations of the first type using the BMA. The relations of type
$x_i-g$ are computed using standard linear algebra techniques. We
have implemented the algorithm, which we call EssBM, in
\emph{Macaulay~2}.

The paper is organized as follows.  First we describe the EssBM
algorithm.  In Section 3, we provide the theoretical support for the
algorithm and include a complexity analysis.  In Section 4, we
compare its performance to the BMA, as implemented in
\emph{Macaulay~2}. We conclude our paper with a discussion of future
directions.

\section{The EssBM Algorithm}

Let $R=k[x_1,\ldots ,x_n]$ where $k$ is a field, and $\se$ be a
fixed term order on $R$. Consider a variety $V\subset k^n$ of points
with multiplicity one and $|V|=m<\infty$.  Here we are primarily
interested in finite fields, where these conditions will
automatically be satisfied for all varieties.  The goal of the EssBM
algorithm is to construct the reduced \gb $G$ with respect to~$\se$
for the ideal $\mathbf{I}(V)$ of points in $V$ and the set
$\mathfrak{B}(G)$ of standard monomials associated to $G$, which
forms a basis for the $k$-vector space $R/\mathbf{I}(V)$. The
algorithm constructs a set $EV\subset \{x_1\ldots ,x_n\}$ of
\emph{essential variables}, a set $SM$ of monomials on $\{x_1\ldots
,x_n\}$, and subsets $GB$ and $Rel$ of the ring $R$. We will see
below that $G$ will be given by $GB\cup Rel$ and $\mathfrak{B}(G)$
by the set $SM$. The \emph{support} (defined in the next section) of
the elements in $SM$ is the set $EV$. We let $EV_i$, $SM_i$, $GB_i$,
and $Rel_i$ denote the $i$-th approximations of the corresponding
sets.

Initialize each set as follows: $EV_0=\{\}$, $SM_0=\{1_R\}$,
$GB_0=\{\}$, and $Rel_0=\{\}$. Let $[n]$ denote the set $\{1,\ldots
,n\}$ and $x^a$ the monomial $x_1^{a_1}\cdots x_n^{a_n}$.  For each
$i\in [n]$, do the following. Find the $i$-th smallest variable, say
$x_i$. Suppose there are $r$ monomials $x^{a_1},\ldots,x^{a_r}$ in
$SM_{i-1}$ that are smaller than $x_i$ in the given ordering.  Try
to write $x_i$ as a $k$-linear combination of these monomials. That
is, find (if they exist) $c_1,\ldots ,c_r\in k$, where
\begin{equation}
\label{linsys}
\begin{split}
  x_i(1) &= \sum_{j=1}^r c_jx^{a_j}(1) \\
  x_i(2) &= \sum_{j=1}^r c_jx^{a_j}(2) \\
  &\cdots  \\
  x_i(m) &= \sum_{j=1}^r c_jx^{a_j}(m)
\end{split}
\end{equation}
and $x^a(t)$ is the evaluation of $x^a$ at the $t$-th point in $V$
for $t\in[m]$. If there are such coefficients, then
$$x_i(t) - \sum_{j=1}^r c_jx^{a_j}(t) = 0$$
for every $t\in[m]$ and it follows that $h:=x_i - \sum_{j=1}^r
c_jx^{a_j}\in \mathbf{I}(V)\cap k[EV_{i-1}\cup \{x_i\}]$, where
$k[EV_{i-1}\cup \{x_i\}]$ is the polynomial ring in the variables in
$EV_{i-1}\cup \{x_i\}$. Since the monomials $x^{a_j}$ were chosen so
that $x_i\se x^{a_j}$, it follows that $x_i$ is the leading term of
an element of $\mathbf{I}(V)$ and so is not a standard monomial.  In
this case let $Rel_i = Rel_{i-1}\cup \{h\}$. If there is no solution
to the system in (\ref{linsys}), then $x_i$ is a standard monomial.
In this case let $EV_i=EV_{i-1}\cup \{x_i\}$, and compute the \gb
$GB_{i}$ and the set $SM_{i}$ of standard monomials for the ideal
$\mathbf{I}(V)\cap k[EV_{i}]$ of the points projected onto the
variables in $EV_{i}$. When $i=n$, return the sets $G:=GB_n\cup
Rel_n$ and $\mathfrak{B}(G):=SM_n$.

Below we give pseudo-code for the complete algorithm, which has been
implemented in \emph{Macaulay~2}. While the BMA computes
\emph{separators} for the points in $V$ in addition to the \gb and
the set of standard monomials, the implementation in
\emph{Macaulay~2} does not. In order to appropriately compare the
two implementations, we do not include separators in this version of
EssBM. However, our algorithm can easily be modified to return the
separators at an additional cost of $O(m)$.

For simplicity, let $[x_j(t)]_{t=1}^m$ denote the $(m\times
1)$-column vector
$$\left(\begin{array}{c}
          x_j(1) \\
          x_j(2) \\
          \vdots \\
          x_j(m)
        \end{array}
\right).$$

\medskip

\noindent \rule{450pt}{0.5pt} \bigskip

\noindent \hspace{6pt}\textbf{The EssBM Algorithm}

\begin{tabbing}
\hspace{8pt}\=\textbf{Input}: \hspace{9pt}\= $V$ a variety;
\hspace{5pt}$\se$ a term order\\

\> \textbf{Output}: \>$G$ the reduced \gb for $\mathbf{I}(V)$ with
respect to $\se$; \\ \> \> $\mathfrak{B}(G)$ the set of standard
monomials for $G$
\end{tabbing}
\noindent \rule{450pt}{0.5pt}

\begin{enumerate}
    \item Initialize: $EV_0:=\{\}; \hspace{5pt} SM_0:=\{1_R\}; \hspace{5pt} GB_0:=\{\}; \hspace{5pt}
    Rel_0:=\{\}$.
    \item For $i$ from 1 to $n$ do
    \item \hspace{12pt} $x_i := i$-th smallest variable
    \item \hspace{12pt} $S := k[EV_{i -1}\cup \{x_i\}]$ with term
    order $\se_S$ induced by $\se$
    \item \hspace{12pt} $r:=|SM_{i-1}|$ and $LM_i:=\{x^{a_j}\pe_S x_i : x^{a_j}\in SM_{i -1}, 1\leq j \leq r\}$ the standard monomials less than~$x_i$
    \item
        \begin{tabbing}
            \hspace{12pt}
            \=$A_i :=(m\times (s+1))$-matrix with $s=|EV_{i -1}|$ \\
            \>first column $[x_i(t)]_{t=1}^m$ \hspace{5pt} and \hspace{5pt} $s$
            columns $[x_j(t)]_{t=1}^m$ for all $x_j\in EV_{i -1}$
        \end{tabbing}
    \item \hspace{12pt} $Eval_i:=(m\times r)$-matrix $(x^{a_j}(p_t))$, where $x^{a_j}\in LM_i$ is evaluated on $p_t$, the point in row $t$ of $A_i$
    \item \hspace{12pt} If there is a solution $c=(c_1,\ldots ,c_r)^T$ to the system of linear equations $Eval_i \cdot c=[x_i(t)]_{t=1}^m$
    \item \hspace{22pt} then $Rel_{i}:=Rel_{i -1}\cup \{x_i-\sum c_jx^{a_j}\}$ where
    $x^{a_j}\in LM_i$
    \item \hspace{22pt} else $EV_{i}:=EV_{i -1}\cup \{x_i\}$ and compute $GB_i$ and $SM_i$ in $k[EV_{i}]$ using the BMA on
    $A_i$
    \item Return $G=GB_n\cup Rel_n$ and $\mathfrak{B}(G)=SM_n$
\end{enumerate}

\medskip

The variables in $EV_n$ are called \emph{essential}. The polynomial
$x_i-\sum c_jx^{a_j}$ computed in the $i$-th step of the algorithm
has $x_i$ as its leading term since the monomials $x^{a_j}$ were
chosen to be smaller than $x_i$. The variables $x_i$ are called
\emph{inessential} since they can be written in terms of essential
variables.

\section{Theoretical Background}

In this section, we provide a detailed proof of the correctness and
worst-case time complexity of the EssBM algorithm. Before stating
and proving the main results, namely Theorems \ref{correctness-gb},
\ref{correctness-sm} and \ref{complexity}, we begin with some
preliminaries.

Recall that the matrix $A_i$ has rows corresponding to the points in
$V$ projected onto the coordinates defined by $EV_i=EV_{i-1}\cup
\{x_i\}$.  Let $P_i$ be this set of projected points.

For the remainder of this paper, we use the shorthand notation $I$
for the ideal $\mathbf{I}(V)$ and $k[EV_i]$ for the polynomial ring
in the variables in the set $EV_i$.  Also, we let $G=GB_n\cup Rel_n$
and $\mathfrak{B}(G)$ the set of standard monomials for $G$.

\begin{lemma}
    The equality $\mathbf{I}(P_i)=I\cap k[EV_i]$ holds.
\end{lemma}

\begin{proof}
    This follows immediately from the construction of the ideal $\mathbf{I}(P_i)$.
\end{proof}

\begin{corollary}
\label{gb}
    The set $GB_i$ is the reduced \gb for the ideal $I\cap k[EV_i]$ with respect to $\se$
    and $SM_i$ is the set of standard
    monomials for $I\cap k[EV_i]$ with respect to $GB_i$.  In
    particular, the statement holds for $i=n$.
\end{corollary}

\begin{proof}
    The sets $GB_i$ and $SM_i$ are the reduced \gb and the set of standard
    monomials, respectively, for the ideal $\mathbf{I}(P_i)$ in $k[EV_i]$.
    From the previous lemma, we have that $\mathbf{I}(P_i)=I\cap k[EV_i]$.
    Hence the result follows.
\end{proof}

Let $f\in R$ be a polynomial.  We define the \emph{support} of $f$,
denoted by $supp(f)$, to be the set of variables that appear in $f$.
By construction, $supp(f)$ is the smallest set $X\subset
\{x_1,\ldots ,x_n\}$ such that $f\in k[X]$. The \emph{support} of a
set of polynomials $S$ is the union over the support of each
polynomial $g\in S$.  Let $LT(f)$ denote the leading term of $f$
with respect to a given term order. The \emph{tail} of $f$ is the
polynomial $tail(f):=f-LT(f)$.

\begin{lemma}
\label{iness}
    Let $f\in R$ be such that $supp(f)\subset EV_n\cup \{x_{\beta_1},\ldots ,x_{\beta_s}\}$
    where $x_{\beta_1}\prec \cdots \prec x_{\beta_s}$ are inessential variables.
    Suppose that  $supp(LT(f))\subset EV_n$.  Then there is
    $f^*\in R$ such that $supp(f^*)\subset EV_n\cup \{x_{\beta_1},\ldots ,x_{\beta_{s-1}}\}$,
    the polynomial $f^*$ has the same leading term as $f$, and $f-f^* \in I$.
\end{lemma}

\begin{proof}
    Consider the largest inessential variable $x_{\beta_s}$. We can
    write
    $$f=LT(f)+\sum_{i=0}^r(x_{\beta_s})^ih_i$$
    where $supp(h_i)\subset EV_n\cup \{x_{\beta_1},\ldots ,x_{\beta_{s-1}}\}$.
    As $x_{\beta_s}$ is an inessential variable, there is an element
    $x_{\beta_s}+g $ of $Rel_n$ with leading term $x_{\beta_s}$.
    Note that $supp(g)\subset EV_n$. Define the polynomial $f'$ from $f$
    by replacing each $(x_{\beta_s})^i$ with~$-(x_{\beta_s})^{i-1}g$:
    $$f'=LT(f)-\sum_{i=0}^r(x_{\beta_s})^{i-1}g h_i.$$
    Then
    $$f-f' = \sum_{i=0}^r\left((x_{\beta_s})^i+(x_{\beta_s})^{i-1}g\right) h_i \in I$$
    since $(x_{\beta_s})^i+(x_{\beta_s})^{i-1}g=(x_{\beta_s})^{i-1}(x_{\beta_s}+g)\in I$.
    As $LT(f)\succ x_{\beta_s} \succ LT(g)$, we have that $LT(f')=LT(f)$. Let~$f^*$
    be the polynomial obtained after $r$ replacements of $x_{\beta_s}$. Note that we
    have $f-f^*\in I$ and $LT(f^*)=LT(f)$.  Since we have replaced all occurrences of $x_{\beta_s}$,
    it follows that $supp(f^*)\subset EV_n\cup \{x_{\beta_1},\ldots
    ,x_{\beta_{s-1}}\}$.
\end{proof}

This lemma gives us a way of removing inessential variables from a
polynomial in $I$ without affecting its leading term, which will be
useful for proving the correctness of EssBM (Theorems
\ref{correctness-gb} and \ref{correctness-sm}). In fact, we can
remove \emph{all} inessential variables.  We emphasize this fact
with the following corollary.

\begin{corollary}
\label{iness2}
    Let $f\in R$. Then there is $f^*\in R$
    such that $supp(f^*)\subset EV_n$, $LT(f^*)=LT(f)$, and $f-f^* \in I$.
\end{corollary}

\begin{theorem}
\label{correctness-gb}
    The set $G$ is the reduced \gb for $I$ with respect to $\se$.
\end{theorem}

\begin{proof}
We first show that $G\subset I$. Consider $g\in G$.  If $g\in GB_n$,
then $g\in I$.  Suppose that $g\in Rel_n$.  Then $g$ is of the form
$x_i - \sum c_jx^{a_j}$ for some $c_j\in k$ and $x^{a_j}\in
R=k[x_1,\ldots ,x_n]$. The coefficients $c_j$ were chosen so that
$x_i(t)=\sum c_jx^{a_j}(t)$ for all $t\in [m]$. Therefore by
construction $g\in I$.

Now let $f\in I$.  We must show that there is some $g\in G$ such
that $LT(g)\mid LT(f)$.  We distinguish two cases.

\begin{itemize}
    \item[Case 1:] $supp(LT(f))\not \subset EV_n$.

    Suppose that $LT(f)$ contains an inessential variable $x_i$.
    By construction of the set $Rel_n$, there is an element~$g$ of
    $Rel_n\subset G$ with leading term $x_i$. It follows that
    $LT(g)$ divides the leading term of~$f$.

    \item[Case 2:] $supp(LT(f))\subset EV_n$.

    Recall that the set $GB_n$ is a \gb of the projection of $I$ onto
    the variables in $EV_n$ (see Corollary~\ref{gb}). If
    $supp(tail(f))$ is also contained in $EV_n$, then $f\in k[EV_n]$
    and there is a $g\in GB_n\subset G$ whose leading term divides $LT(f)$.

    Assume that $supp(tail(f))\not \subset EV_n$.  Using Corollary \ref{iness2},
    we can find $h\in I$ such that $supp(f-h)\subset
    EV_n$ and $LT(f-h)=LT(f)$.  Since $f-h\in k[EV_n]$, there is a
    $g\in GB_n\subset G$ whose leading term divides $LT(f-h)=LT(f)$.
\end{itemize}

To prove that $G$ is reduced, let $g\ne h\in G$.  We wish to show
that $g$ and $h$ satisfy the following criterion:
\begin{equation}
    \label{div}
    LT(g)\text{ does not divide any monomial in }h.
\end{equation}
We consider the following four cases.
\begin{itemize}
    \item[Case 1:] $g,h\in GB_n$.

    As $GB_n$ is the reduced \gb for the ideal $I$ projected onto the essential variables,
    then $g,h$ satisfy~(\ref{div}).

    \item[Case 2:] $g,h\in Rel_n$.

    Let $LT(g)=x_i$ and $h=x_j-\sum_{i}c_{i}x^{a_i}$ for $i\neq j$.  Note
    that $supp(h)\subset EV_{j-1}\cup \{x_j\}$.  Clearly $x_i$ does not
    divide $x_j$.  As $supp(tail(h))\subset EV_n$ and $x_i\notin EV_n$, then
    $x_i$ does not divide any monomials in $tail(h)$.

    \item[Case 3:] $g\in GB_n$ and $h\in Rel_n$.

    Let $LT(h) = x_i$ for some inessential variable. This will not be
    divisible by $LT(g)$, which contains at least one essential
    variable.  All other terms $x^a$ of $h$ are standard monomials
    for the projection of $I$ onto the variables in $EV_i$; in
    particular, $supp(x^a) \subset EV_i$.  It follows that if
    $supp(g) \subseteq EV_i$, then $LT(g)$ does not divide any term of
    $h$.  By Corollary \ref{gb}, $supp(g)$ contains only essential variables.  Thus if
    $supp(g)$ is not contained in $EV_i$, then $supp(g)$ must contain a
    variable $x_j$ with $x_i \prec x_j$.  This $x_j$ divides some term
    $x^b$ of $g$, and it follows that if $LT(g)$ divides some term
    $x^a$ of $h$, then $x_j \preceq x^b \preceq LM(g) \preceq
    x^a \preceq x_i$, which contradicts the assumption that $x_i
    \prec x_j$.

    \item[Case 4:] $g\in Rel_n$ and $h\in GB_n$.

    Then $LT(g)$ is some inessential variable, say $x_i$.  However, $supp(h)\subset EV_n$
    and so $g,h$ satisfy criterion (\ref{div}).
\end{itemize}

\end{proof}

Next we compute the number of elements in $\mathfrak{B}(G)$ and show
the relationship between $\mathfrak{B}(G)$ and the set $SM_n$.

\begin{lemma}
\label{bv}
    The set $\mathfrak{B}(G)$ has $|V|$ elements.
\end{lemma}

The previous lemma is usually stated for algebraically closed fields
$k$ and proved with the help of the Strong Hilbert Nullstellensatz
(see \cite{cox2}). We include a proof of the statement for the case
where all points have multiplicity one, as is being assumed
throughout the paper.

\begin{proof}
    Suppose $V=\{a_1,\ldots ,a_m\}$ and define $I_i:=\mathbf{I}(\{a_i\})$.
    Then
    $I=\mathbf{I}(\bigcup_{i=1}^m \{a_i\})=\bigcap_{i=1}^m I_i$, since each
    point $a_i$ has multiplicity one. Note
    that each of the ideals $I_i$ is maximal and it follows that they are pairwise comaximal.
    Consider the quotient ring $R/I$.  By the Chinese Remainder Theorem,
    there is a ring homomorphism such that
    $$R/I\cong R/I_1 \times \cdots \times R/I_m.$$
    As each $I_i$ is maximal, then each $R/I_i\cong
    k$ and it follows that $R/I\cong k^m$, as rings. Further, the quotient ring and~$k^m$
    can be viewed as $k$-vector spaces, and the isomorphism can be extended to an isomorphism of vector
    spaces.  Hence, the dimension of $R/I$ as a vector space is $dim_k(R/I)=m$.
    Since $\mathfrak{B}(G)$ forms a basis for the vector space $R/I$ (Proposition 2.1.6 in \cite{adams}),
    we conclude that $|\mathfrak{B}(G)|=m=|V|$.
\end{proof}

\begin{theorem}
\label{correctness-sm}
    The set $SM_n$ is the set of standard monomials for $I$ with respect to $G$.
\end{theorem}

\begin{proof}
    By Corollary \ref{gb}, we have that $SM_n$ is the set of standard monomials
    for the ideal $I\cap k[EV_n]$ with respect to the \gb $GB_n$.  As $V$ has
    finitely many points, then $|\mathfrak{B}(G)|=|V|$.  Consider
    a monomial $x^a\in \mathfrak{B}(G)$.  If $x^a\not \in k[EV_n]$,
    then it contains an inessential variable, say $x_i$.  As $x_i$
    is the leading term of an element in $Rel_n\subset I$, it is not
    a standard monomial for $G$, contradicting the assumption that $x^a\not \in
    \mathfrak{B}(G)$.  Therefore $x^a \in k[EV_n]$.

    By construction, $x^a\not \in LT(I)$.  Using the
    set-containment relation
    $$LT\left(I\cap k[EV_n]\right)\subset LT(I),$$
    it follows that $x^a\not \in LT(I\cap k[EV_n])$ and so $\mathfrak{B}(G)\subset SM_n$.
    To see equality, note that the set $P_n$ of projected points defined by $EV_n$ has at
    most as many points as $V$.  Then $|SM_n|=|P_n|\leq |V|=m$.
    Since $\mathfrak{B}(G)\subset SM_n$, it follows that $m=|\mathfrak{B}(G)|\leq |SM_n|\leq m$.
    Hence $\mathfrak{B}(G)=SM_n$; that is, $SM_n$ is the set of standard monomials for $I$ with respect to $G$.
\end{proof}

We conclude this section with a complexity analysis of EssBM.

\begin{theorem}
\label{complexity}
    The EssBM algorithm terminates and has worst-time
    complexity $O(nm^3) + O(m^6)$, which is dominated by $O(nm^3)$ when
    $m\ll n$.
\end{theorem}

\begin{proof}
    We compute the complexity of each step and then provide a
    summary at the end. Step~1 has complexity $O(1)$.
    In Step~2, the algorithm enters a loop of length $n$.  Steps~3-8 are executed
    in each iteration of the loop.  They have the following complexities:
    \begin{itemize}
        \item[Step 3.] $O(1)$: Executing this step requires constant time
        since the variable order, given as part of the declaration of the
        term order, is maintained in one array.
        \item[Step 4.] $O(m^2)$: This step may not even be required by all implementations;
        if required, it involves passing $O(m^2)$ variables to a new object of size $O(m^2)$.
        \item[Step 5.] $O(m^3)$: As term orders are typically stored
        as matrices, in this case the term order $\pe_S$ is a matrix of
        dimension $O(m^2)$.  Determining the order between two monomials of $S$
        requires multiplication of a vector of length $O(m)$ by this matrix.
        So for each monomial $x^a\in SM_{i-1}$, there are at most $m^2$ operations
        required for comparing $x_i$ to $x^a$ and there are at most $m$ such monomials.
        \item[Step 6.] $O(m)$: An $(m\times 1)$-column vector is added to a matrix with columns
        corresponding to the variables in $EV_{i-1}$.
        \item[Step 7.] $O(m^3)$: As there are at most $m$ variables
        in each monomial and at most $m^2$ entries in the matrix, the cost of
        executing this step is $O(m^3)$.
        \item[Step 8.] $O(m^3)$: Solving a linear system of $m$ equations in
        $r\leq m$ unknowns requires $O(m^3)$ time.
      \end{itemize}

        Step~9 has complexity $O(1)$ and will be executed at most $n$ times.

        Since there can be at most $m$ essential variables, Step~10 will be executed
        at most~$m$ times. The complexity of each execution of Step~10 is $O(m^5)$:
        Updating $EV_{i}$ is a constant operation.  However, computing $GB_i$ and $SM_i$
        for the matrix $A_{i}$ is associated to the cost of
        calling the BMA, which is quadratic in the number
        of variables and cubic in the number of points.  In this case,
        the numbers of variables and points are given by the dimensions
        of~$A_{i}$.  Since both row and column dimensions are bounded
        above by~$m$, it follows that the complexity of executing this
        step is $O(m^5)$.

        Step~11 has complexity $O(n + m^2)$: Note that there are $O(m^2)$ elements in $GB_n$
        (see \cite{mariani}), $O(n)$ relations in
        $Rel_n$, and $m$ monomials in $SM_n$. So returning these
        sets requires $O(n + m^2 + m)$ operations.

    Hence, we can calculate the total complexity $C(EssBM)$ of the algorithm as follows:
    \begin{eqnarray}
    \nonumber
        C(EssBM) &=& O(1)+O(n)\left[ O(1 + m^2 + m^3+ m + m^3 + m^3 + 1)\right] + O(m)O(m^5) + O(n + m^2 + m) \\
    \nonumber
        &=& O(nm^3) + O(m^6).
    \end{eqnarray}
      When $m\ll n$, then $O(nm^3)$ becomes the dominating term and the
    above estimate reduces to
    $$C(EssBM)=O(nm^3).$$
\end{proof}

\section{Performance of the EssBM Algorithm}

To test the performance of our algorithm, we compared its run-time
to that of the BMA\footnote{The Buchberger-M\"oller algorithm has
been implemented as the function \emph{points} in the ``Points''
package of \emph{Macaulay 2} distribution version 0.9.8.}, as
implemented in \emph{Macaulay 2}, on randomly generated varieties in
$k^n$. For this analysis, we let the field $k$ be $\mathbb{F}_p$ for
$p\in \{3,17\}$. Since the complexities of the two algorithms depend
on $m$ and $n$, we chose a range of values for these parameters,
namely, $m\in \{5,10,15\}$ and $n\in \{100,150,200,250,300\}$.  For
each set of parameters $p,m,$ and $n$, we generated 10 varieties
using a built-in random number generator in \emph{Macaulay 2},
without specifying prior constraints on the relative position of the
points in the variety. We performed this experiment using two term
orders: a lexicographic order (\emph{lex}) and a graded reverse
lexicographic order (\emph{grevlex}), each with the same variable
order.

\begin{figure}[h]
  \centering
  \caption{Run-times averaged over 10 randomly generated varieties for $p=3$ and \emph{lex}.}
  \label{lex}
  \includegraphics[width=3.5in,angle=270]{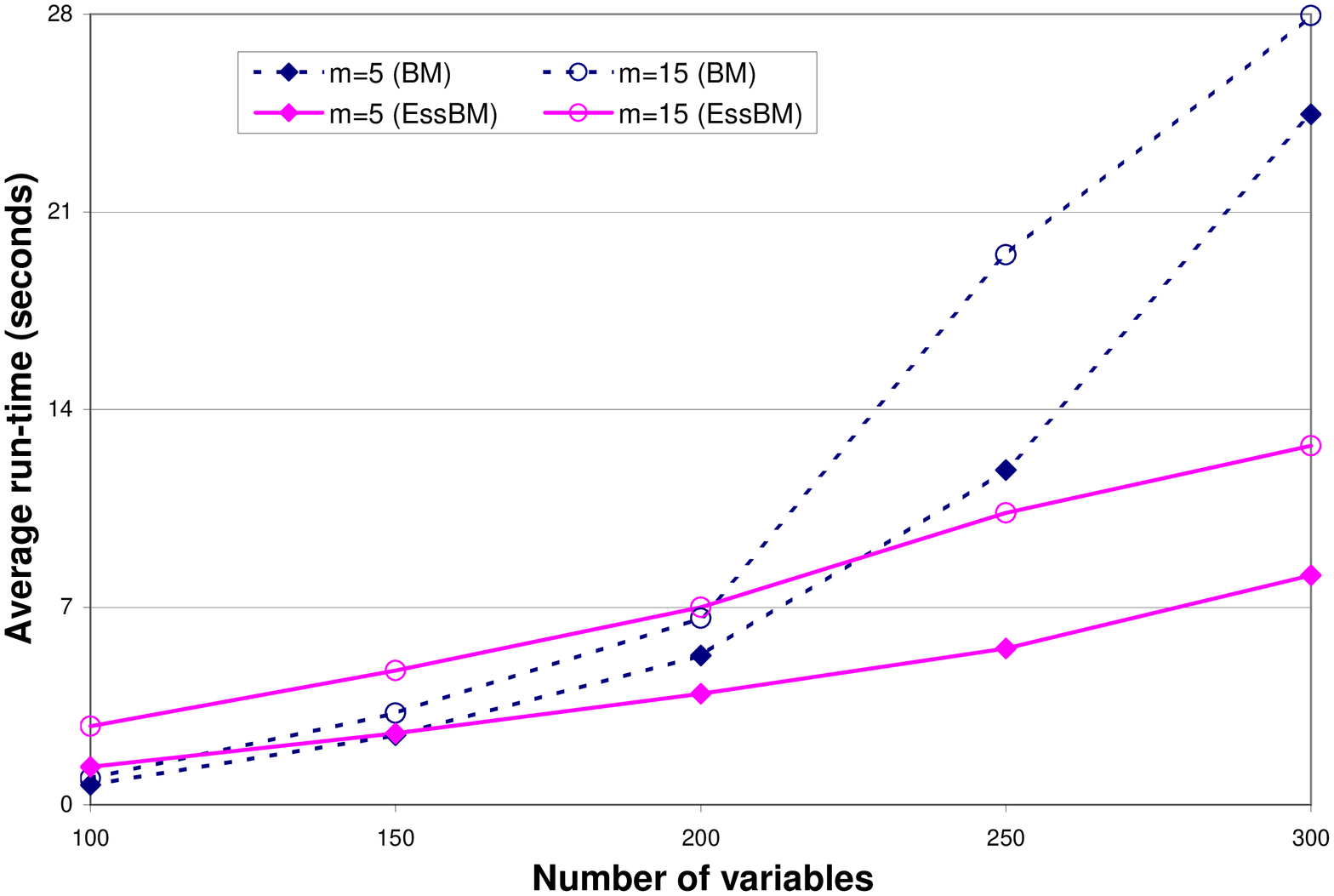}
\end{figure}

\begin{figure}[b!]
  \centering
  \caption{Run-times averaged over 10 randomly generated varieties for $p=3$ and \emph{grevlex}.}
  \label{grevlex}
  \includegraphics[width=3.5in,angle=270]{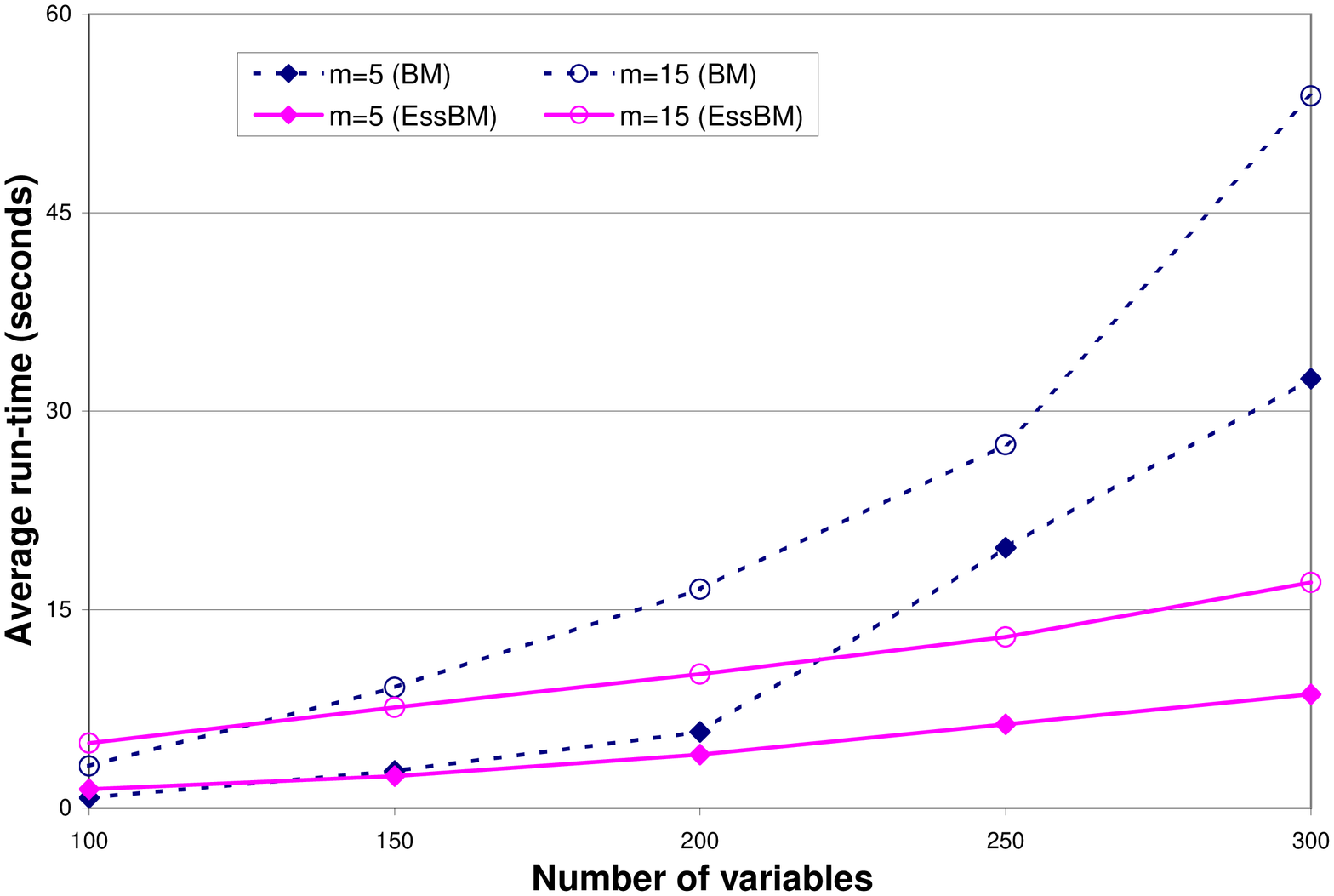}
\end{figure}

Figures \ref{lex} and \ref{grevlex} show the run-times for the two
algorithms for $p=3$ and $m=5,15$.  As the run-times for $m=10$ fall
between the $m=5$ and $m=15$ settings, we omitted them from the
plots.  We display the results for all parameters settings in the
appendix.  The run-times for $p=17$ are similar.

As a measure of the stability of the run-time data, we computed the
\emph{coefficient of variation}, defined to be the ratio of the
standard deviation to the mean of the data. For the \emph{grevlex}
experiments, this coefficient ranges from 0.004 to 0.2, whereas for
the \emph{lex} experiments it ranges from 0.01 to 0.1. Since this
implies very low variability of the run-times for fixed $p,n,$ and
$m$, we displayed only mean values in Figures \ref{lex} and
\ref{grevlex}.

The empirical results corroborate our theoretical prediction that
for $m\ll n$, the EssBM algorithm outperforms the BMA.  For small
$n$, however, we observe that EssBM is slower, which we attribute to
the overhead costs associated to multiple calls to the BMA.

\section{Discussion}

Recently, applications of \gbs as a promising model selection tool
in molecular biology have been proposed \cite{allen,LS}. These
applications require computation of a \gb for a zero-dimensional
ideal $\mathbf{I}(V)$ in a polynomial ring $k[x_1,\ldots ,x_n]$,
where $|V|=m\ll n$.  Previously, no algorithms for computing \gbs
optimized for $m\ll n$ had been available.  The run-time of the
existing implementations was a bottleneck in applications of the
methods in \cite{allen} and \cite{LS} to data sets whose size is of
the order typical for biochemical data sets such as microarray data.

The EssBM algorithm presented here goes some way towards alleviating
this problem in that it reduces the worst-case complexity, which is
$O(n^2m^3)$ for the standard Buchberger-M\"oller algorithm, to
$O(nm^3)$ for $m\ll n$. Our implementation and testing indicate that
for a small number of distinct points in general position, EssBM
starts outperforming a standard implementation of the BMA when the
number of variables exceeds 200. This should make it possible to use
the methods of \cite{allen} and \cite{LS} for analysis of larger
data sets than was hitherto possible. Unfortunately, the worst-time
complexity estimate $O(nm^3+m^6)$ of the EssBM algorithm suggests
that it may still be infeasible for moderately large $m$. We are
currently working on a related algorithm that would further reduce
this complexity.

\section{Acknowledgements}
We wish to thank Luis Garc{\'i}a-Puente and an anonymous reviewer
for their valuable comments.

This research was done during the visit of WJ to the Mathematical
Biosciences Institute in the academic year 2006/2007 and supported
by the National Science Foundation under Agreement No. 0112050.

\bibliographystyle{amsplain}
\bibliography{mod-bm}

\pagebreak

\appendix
{\Large \textbf{Appendix}}

\begin{figure}[h]
  \centering
  \caption{Run-times for 10 randomly generated varieties and \emph{lex}.}
  \medskip
  \includegraphics[width=6.5in]{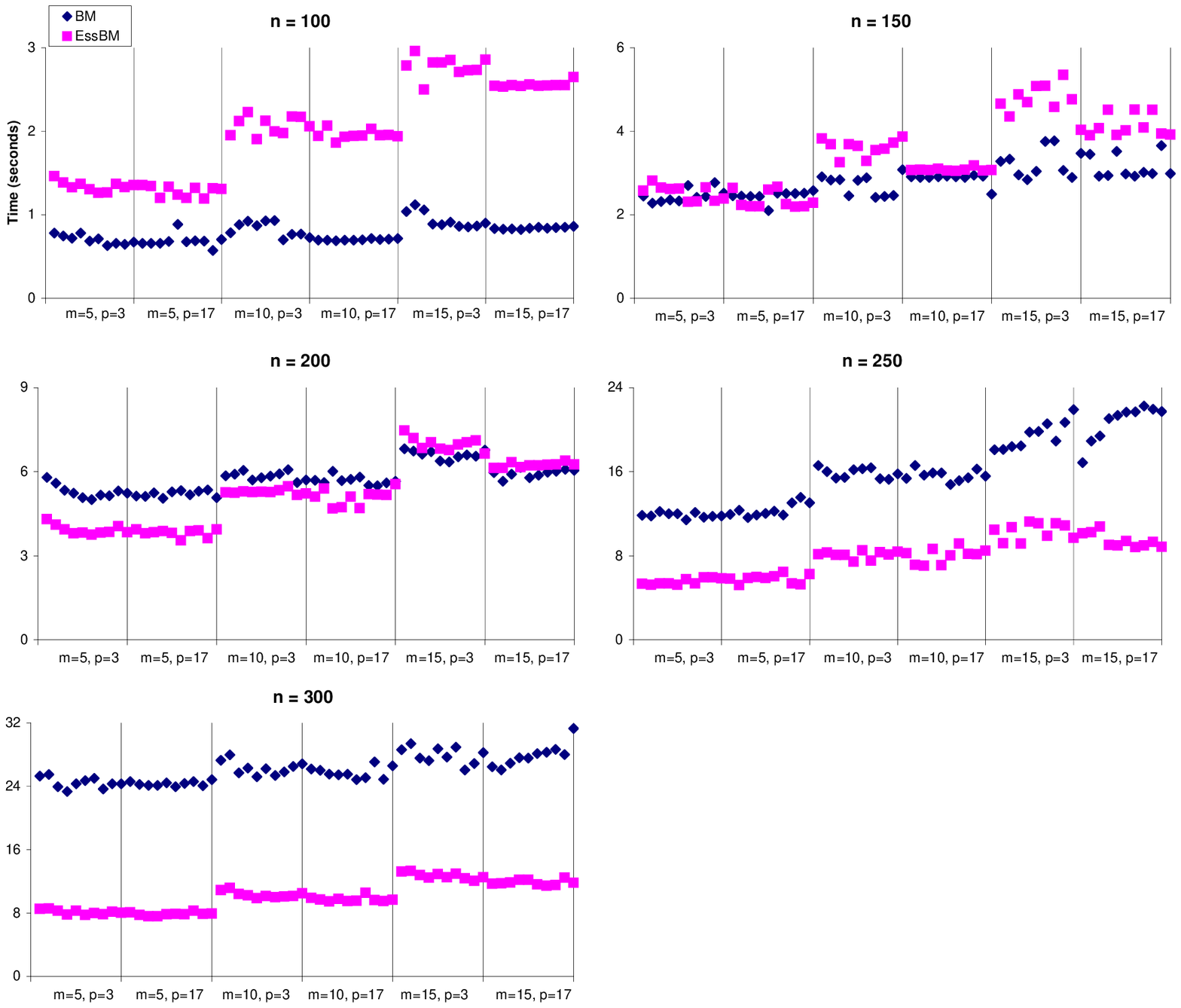}
\end{figure}

\begin{figure}[h]
  \centering
  \caption{Run-times for 10 randomly generated varieties and \emph{grevlex}.}
  \medskip
  \includegraphics[width=6.5in]{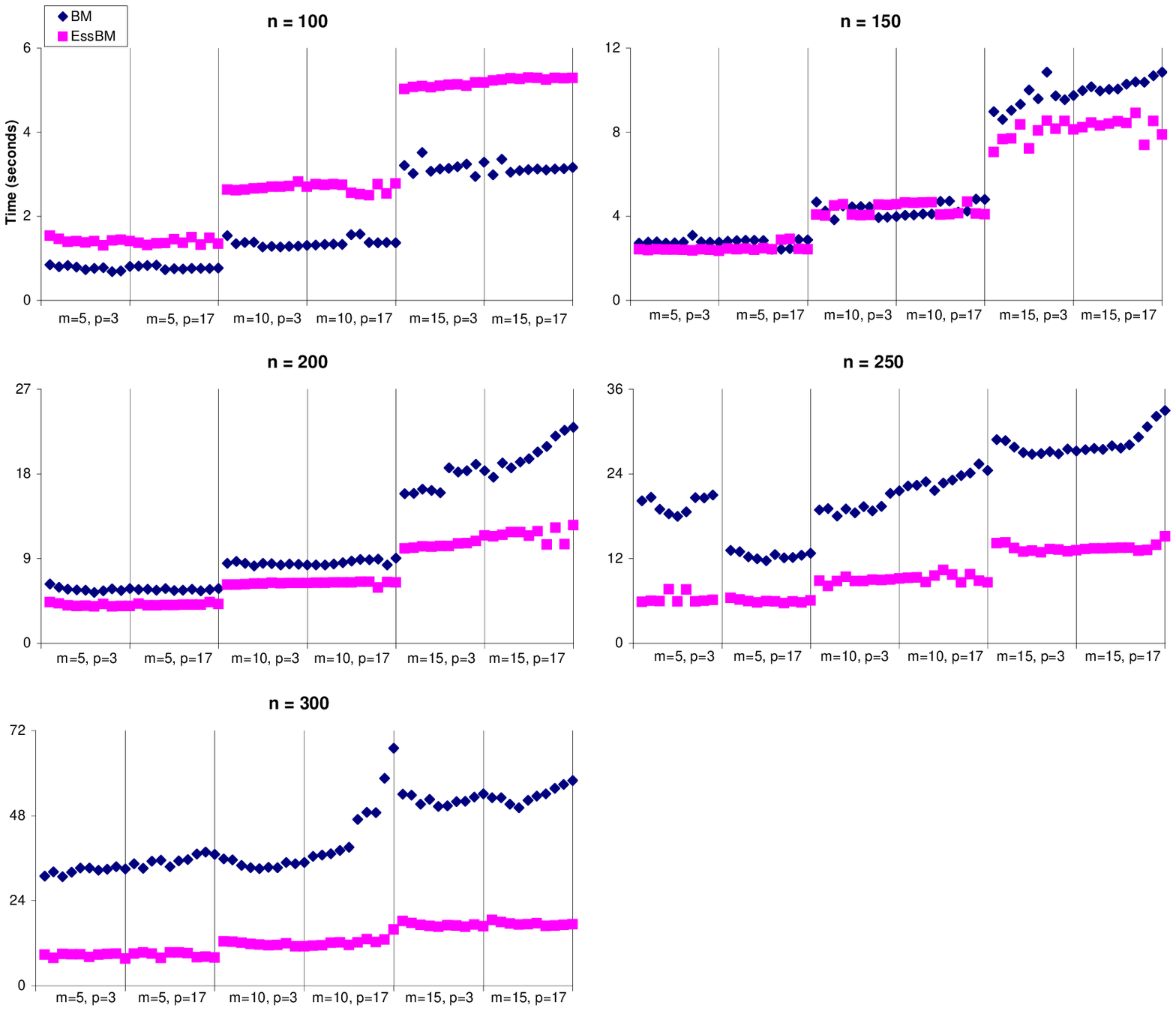}
\end{figure}

\end{document}